\documentclass{article}%
\usepackage{amssymb}
\usepackage{amsmath}
\usepackage{amsfonts}
\usepackage{graphicx}%
\providecommand{\U}[1]{\protect\rule{.1in}{.1in}}
\newtheorem{theorem}{Theorem}

\newtheorem{definition}[theorem]{Definition}

\newtheorem{lemma}[theorem]{Lemma}

\newtheorem{proposition}[theorem]{Proposition}

\begin{document}

\title{Kesten's theorem for Invariant Random Subgroups}
\author{Mikl\'{o}s Ab\'{e}rt, Yair Glasner and B\'{a}lint Vir\'{a}g}
\maketitle

\begin{abstract}
An invariant random subgroup of the countable group $\Gamma$ is a random
subgroup of $\Gamma$ whose distribution is invariant under conjugation by all
elements of $\Gamma$.

We prove that for a nonamenable invariant random subgroup $H$, the spectral
radius of every finitely supported random walk on $\Gamma$ is strictly less
than the spectral radius of the corresponding random walk on $\Gamma/H$. This
generalizes a result of Kesten who proved this for normal subgroups.

As a byproduct, we show that for a Cayley graph $G$ of a linear group with no
amenable normal subgroups, any sequence of finite quotients of $G$ that
spectrally approximates $G$ converges to $G$ in Benjamini-Schramm convergence.
In particular, this implies that infinite sequences of finite $d $-regular
Ramanujan Schreier graphs have essentially large girth.

\end{abstract}

\section{Introduction}

For a $d$-regular, countable, connected undirected graph $G$, we define the
\textit{spectral radius} of $G$, denoted $\rho(G)$, to be the norm of the
Markov averaging operator on $\ell^{2}(G)$. This norm can also be expressed
as
\begin{equation}
\rho(G)=\lim_{n\rightarrow\infty}\left(  p_{x,x,2n}\right)  ^{1/2n}\tag{Rho}%
\end{equation}
where $p_{x,y,k}$ denotes the probability that a simple random walk of length
$k$ starting at $x$ ends at $y$.

Let $\Gamma$ be a group generated by the symmetric set $S$ and let
$\mathrm{Cay}(\Gamma,S)$ denote the Cayley graph of $\Gamma$ with respect to
$S$. Let $H$ be a subgroup of $\Gamma$ and let $\mathrm{Sch}(\Gamma/H,S)$ be
the Schreier graph of $\Gamma$ defined on the coset space $\Gamma/H$. Since
the map $g\mapsto Hg$ extends to a covering map from $\mathrm{Cay}(\Gamma,S)$
to $\mathrm{Sch}(\Gamma/H,S)$, using (Rho) we get
\[
\rho(\mathrm{Sch}(\Gamma/H,S))\geq\rho(\mathrm{Cay}(\Gamma,S))\text{. }%
\]

In his seminal papers \cite{kesten1} and \cite{kesten2}, Kesten proved the
following result.

\begin{theorem}
[Kesten]\label{kestenoriginal}Let $\Gamma$ be a group generated by a finite
symmetric set $S$ and let $N$ be a normal subgroup of $\Gamma$. Then the
following are equivalent: \newline1) $\rho(\mathrm{Cay}(\Gamma/N,S))=\rho
(\mathrm{Cay}(\Gamma,S))$;\newline2) $N$ is amenable.
\end{theorem}

Two special cases of this theorem are particularly well known in the
literature. First, applying the result to $N=\Gamma$ we get that $\Gamma$ is
amenable if and only if it has spectral radius $1$. Second, applying it to the
free group and noting that free groups have no nontrivial amenable normal
subgroups, gives that if a $d$-regular infinite Cayley graph has the same
spectral radius as the $d$-regular tree, then it is isomorphic to it. This has
been later generalized to vertex transitive graphs as well (see e.g.
\cite{paschke}).

It is natural to ask whether Theorem \ref{kestenoriginal} holds without the
normality assumption, replacing $\mathrm{Cay}(\Gamma/N,S)$ with the Schreier
graph $\mathrm{Sch}(\Gamma/H,S)$. It can be shown that if $H$ is amenable,
then $\rho(\mathrm{Sch}(\Gamma/H,S))=\rho(\mathrm{Cay}(\Gamma,S))$. Already
Kesten's original proof does not use normality at this point. However, the
reverse implication (and hence Theorem \ref{kestenoriginal}) fails miserably
for subgroups in general. For instance, when $\Gamma$ is a free group on a
large enough finite set $S$ and $H$ is the subgroup generated by the first $2$
generators of $\Gamma$, then $\mathrm{Sch}(\Gamma/H,S)$ has the same spectral
radius as $\mathrm{Cay}(\Gamma,S)$, although $H$ is obviously not amenable.

Nevertheless, Theorem \ref{kestenoriginal} does hold for a natural stochastic
generalization of normal subgroups.

\begin{definition}
Let $\Gamma$ be a countable group. An \emph{invariant random subgroup} (IRS)
of $\Gamma$ is a random subgroup of $\Gamma$ whose distribution is invariant
under conjugation by $\Gamma$.
\end{definition}

The main result of this paper is the following.

\begin{theorem}
\label{kestenaction}Let $\Gamma$ be a group generated by a finite symmetric
set $S$ and let $H$ be an invariant random subgroup of $\Gamma$. Then the
following are equivalent: \newline1) $\rho(\mathrm{Sch}(\Gamma/H,S))=\rho
(\mathrm{Cay}(\Gamma,S))$ a.s.\newline2) $H$ is amenable a.s.
\end{theorem}

Applying Theorem \ref{kestenaction} to the Dirac measure on a fixed normal
subgroup of $\Gamma$, we get back Theorem \ref{kestenoriginal}.

For a finite, regular graph $G$ we have $\rho(G)=1$. Let $\rho_{0}(G)$ denote
the norm of the Markov averaging operator acting on $\ell_{0}^{2}(G)$, the
space of vectors with zero sum. Let $\Gamma$ be a group generated by a finite
symmetric set $S$ and let $H\leq\Gamma$ be a subgroup of finite index. A
priori, it may happen that $\rho_{0}(\mathrm{Sch}(\Gamma/H,S))<\rho
(\mathrm{Cay}(\Gamma,S))$, but the following analogue of the Alon-Boppana
theorem \cite{alon} shows that asymptotically, $\mathrm{Sch}(\Gamma/H,S)$ can
not beat $\mathrm{Cay}(\Gamma,S)$ spectrally.

\begin{proposition}
\label{abo}Let $\Gamma$ be an infinite group generated by the finite symmetric
set $S$ and let $H_{n}\leq\Gamma$ be a sequence of subgroups of finite index
with $\left\vert \Gamma:H_{n}\right\vert \rightarrow\infty$. Then
\[
\lim\inf\rho_{0}(\mathrm{Sch}(\Gamma/H_{n},S))\geq\rho(\mathrm{Cay}%
(\Gamma,S))\text{.}%
\]

\end{proposition}

We say that the sequence $H_{n}$ of subgroups of finite index \emph{locally
approximates} $\Gamma$, if $\mathrm{Sch}(\Gamma/H_{n},S)$ converges to
$\mathrm{Cay}(\Gamma,S)$ in Benjamini-Schramm convergence, that is, for every
$R>0$ we have
\[
\lim_{n\rightarrow\infty}\frac{\left\vert \left\{  v\in\mathrm{Sch}%
(\Gamma/H_{n},S)\mid B_{R}(v)\cong B_{R}\right\}  \right\vert }{\left\vert
\Gamma:H_{n}\right\vert }%
\]
where $B_{R}(v)$ denotes the ball of radius $R$ centered at $v$ and $B_{R}$
denotes the ball of radius $R$ in $\mathrm{Cay}(\Gamma,S)$. This is equivalent
to the following condition: for every $1\neq g\in\Gamma$, we have
\[
\lim_{n\rightarrow\infty}\frac{\left\vert \mathrm{Fix}(\Gamma/H_{n}%
,g)\right\vert }{\left\vert \Gamma:H_{n}\right\vert }=0
\]
where
\[
\mathrm{Fix}(\Gamma/H,g)=\left\{  Hx\mid Hxg=Hx\right\}
\]
is the set of $H$-cosets fixed by $g$. This in particular shows that the
notion of local approximation is independent of $S$.

It is easy to see that every normal chain with trivial intersection locally
approximates its ambient group. Another family of examples is when $\Gamma$ is
an arithmetic group in zero characteristic and $H_{n}$ is an arbitrary
sequence of congruence subgroups in $\Gamma$. See \cite{abgv} for a proof
using Strong Approximation and \cite{samurai} for explicit estimates on the
essential injectivity radius.

Using \cite{abgv}, Theorem \ref{kestenaction} now implies that for a finitely
generated linear group $\Gamma$ with no nontrivial amenable normal subgroups,
every sequence that spectrally approximates $\Gamma$ must locally approximate it.

\begin{theorem}
\label{farber}Let $\Gamma$ be a finitely generated linear group with no
nontrivial amenable normal subgroups and let $S$ be a finite symmetric
generating set for $\Gamma$. Let $H_{1},H_{2},\ldots$ be a sequence of
subgroups of finite index with $\left\vert \Gamma:H_{n}\right\vert
\rightarrow\infty$ such that
\[
\lim\limits_{n\rightarrow\infty}\sup\rho_{0}(\mathrm{Sch}(\Gamma/H_{n}%
,S))\leq\rho(\mathrm{Cay}(\Gamma,S))\text{.}%
\]
Then $H_{n}$ locally approximates $\Gamma$.
\end{theorem}

Note that the assumption of having no nontrivial amenable normal subgroups is
not very restrictive, as for any finitely generated linear group $\Gamma$, the
group $\Gamma/A(\Gamma)$ is also linear, where $A(\Gamma)$ denotes the maximal
amenable normal subgroup of $\Gamma$.

We can apply Theorem \ref{farber} to the special case when $\Gamma$ is a
nonamenable free group. We call the finite $d$-regular graph $G$
\emph{Ramanujan}, if $\rho_{0}(G)\leq\rho(T_{d})$ where $T_{d}$ denotes the
$d$-regular tree. Lubotzky, Philips and Sarnak \cite{lps}, Margulis
\cite{margulis} and Morgenstein \cite{morgen} have constructed sequences of
$d$-regular Ramanujan graphs for $d=p^{\alpha}+1$. Also, Friedman
\cite{friedman} showed that random $d$-regular graphs are close to being Ramanujan.

All the known Ramanujan graphs have \emph{essentially large girth}, that is,
for all $L$, we have
\[
\lim_{n\rightarrow\infty}\frac{c_{L}(G_{n})}{\left\vert G_{n}\right\vert }=0
\]
where $c_{L}(G_{n})$ denotes the number of cycles of length $L$ in $G_{n}$.
Theorem \ref{farber} shows that this is not a coincidence.

\begin{theorem}
\label{raman}Let $\Gamma$ be a group generated by a finite symmetric set $S$
and let $H_{n}\leq\Gamma$ be a sequence of subgroups of finite index with
$\left\vert \Gamma:H_{n}\right\vert \rightarrow\infty$ such that the Schreier
graphs $\mathrm{Sch}(\Gamma/H_{n},S)$ are Ramanujan ($n\geq1$). Then $\Gamma$
is a free group freely generated by $S$ and $\mathrm{Sch}(\Gamma/H_{n},S)$ has
essentially large girth.
\end{theorem}

The result follows from the original Kesten theorem in the case when the
$H_{n}$ are normal in $\Gamma$, but for arbitrary subgroups, one needs Theorem
\ref{kestenaction}.

Note that in the forthcoming paper \cite{agv}, which focuses on graph theory
and is not assuming the graphs to be Schreier graphs, we prove more: we give
explicit estimates between the spectral radius and the density of short
cycles. In particular, we show that the essential girth of a Ramanujan graph
$G$ is at least $\log\log\left\vert G\right\vert $. \bigskip

\noindent\textbf{Invariant Random Subgroups: history and background.} Although
we are coining the name IRS in this paper, they have been around in various forms.

A random sample of an IRS $H$ is not just an arbitrary subgroup of $\Gamma$.
This can be visualized by looking at the Schreier graph $\mathrm{Sch}%
(\Gamma/H,S)$. This graph exhibits a `statistical homogeneity' in the sense
that every event that occurs in it locally will appear in it again with some
measurable frequency. This phenomenon is pointed out by Aldous and Lyons in
their starting paper on unimodular random networks \cite{aldouslyons}. It
turns out that for a random subgroup $H$ of $\Gamma$, $\mathrm{Sch}%
(\Gamma/H,S)$ forms a unimodular random network if and only if $H$ is an IRS.
See Section 3 for details.

Another source of invariant random subgroups is probability measure preserving
actions: the stabilizer of a random point of the underlying space always forms
an IRS, and vice versa, every IRS arises this way (see Proposition
\ref{ekviv}). A result of St\"{u}ck and Zimmer \cite{stuckzimmer} says that in
a higher rank simple real Lie group $G$, every ergodic IRS equals a
Haar-random conjugate of a lattice in $G$. This generalizes the Margulis
normal subgroup theorem \cite{margulis} and is the deepest result on invariant
random subgroups so far. Although the result itself is not about countable
groups, it implies that in a \emph{lattice} $\Gamma$ in a Lie group as above,
like $SL_{3}(\mathbb{Z})$, every ergodic IRS is either trivial or has finite
index in $\Gamma$. Since these lattices satisfy the congruence subgroup
property, this means that we have a complete control on their invariant random subgroups.

The stabilizer of a random point of a probability measure preserving action
has also been studied by Bergeron and Gaboriau \cite{berggab}, who pointed out
that these groups tend to behave similarly to normal subgroups. They present
this in \cite[Theorem 5.4]{berggab} where they prove that if $\Gamma$ has
positive first $L^{2}$ Betti number, then any IRS in $\Gamma$ which is
nontrivial and of infinite index has infinite first $L^{2}$ Betti number.
Another recent result, using the language of measure preserving actions is due
to Vershik \cite{vershik}, who classified invariant random subgroups of the
finitary symmetric group of countable rank.

Exploiting the St\"{u}ck-Zimmer result, it is proved in \cite{samurai} that
for a higher rank simple real Lie group with symmetric space $X$, every
sequence of symmetric $X$-manifolds with volume tending to infinity must
converge to $X$ in Benjamini-Schramm convergence. With some additional
assumptions, this implies convergence of normalized Betti numbers and
multiplicities of certain unitary representations. Also, in \cite{abgv} it is
shown that in a linear group $\Gamma$, every amenable IRS of $\Gamma$ lies in
the amenable radical, the maximal amenable normal subgroup of $\Gamma$. Note
that we do not know any counterexamples to this phenomenon, even when we do
not assume linearity. This result is used in Theorem \ref{farber} but not in
Theorem \ref{raman}, as it is easy to see that non-Abelian free groups do not
possess nontrivial amenable invariant random subgroups.

\section{Preliminaries \label{section_preliminaries}}

In this section we define the notions and state some basic results used in the
paper.\bigskip

A subset $S$ of the group $\Gamma$ is called \emph{symmetric}, if for all
$s\in S$ we have $s^{-1}\in S$. Let $\Gamma$ be a group generated by the
finite symmetric subset $S$ and let $H$ be a subgroup of $\Gamma$. We define
the \emph{Schreier graph} $\mathrm{Sch}(\Gamma/H,S)$ as follows: the vertex
set is the right coset space $\left\{  Hg\mid g\in\Gamma\right\}  $ and for
each $s\in S$ and vertex $Hg$, there is an $s$-labeled edge going from $Hg$ to
$Hgs$. This defines a directed graph where the $s$-labeled edges are inverses
of the $s^{-1}$-labeled edges. The \emph{root} of $\mathrm{Sch}(\Gamma/H,S)$
is defined as the trivial coset $H$.

When $H$ is the trivial subgroup of $\Gamma$, the Schreier graph is
particularly nice, and we call it the \emph{Cayley graph} $\mathrm{Cay}%
(\Gamma,S)$ of $\Gamma$ with respect to $S$. Cayley graphs are vertex
transitive and carry geometric information on the group $\Gamma$. The map
$g\longmapsto Hg$ extends to a covering map from $\mathrm{Cay}(\Gamma,S)$ to
$\mathrm{Sch}(\Gamma/H,S)$ and so $\mathrm{Cay}(\Gamma,S)$ can also be thought
of as the `universal cover relative to $\Gamma$' of $\mathrm{Sch}(\Gamma/H,S)$
as $H$ varies over all subgroups of $\Gamma$.

Let $G=\mathrm{Sch}(\Gamma/H,S)$. Let $\ell^{2}(G)$ be the Hilbert space of
all square summable functions on the vertex set of $G$. Let us define the
Markov operator $M:\ell^{2}\rightarrow\ell^{2}$ as follows:
\[
(Mf)(x)=\frac{1}{d}%
{\displaystyle\sum\limits_{s\in S}}
f(xs)
\]

We define the \textit{spectral radius} of $G$, denoted $\rho(G)$, to be the
norm of $M$.

For a graph $G$ and $x,y\in V(G)$ let $P_{x,n}$ denote the set of walks of
length $n$ starting at $x$ and let $P_{x,y,n}$ denote the elements in
$P_{x,n}$ that end at $y$. A \emph{random walk of length }$n$\emph{\ starting
at }$x$ is a uniform random element of $P_{x,n}$. Let the \emph{probability of
return}
\[
p_{x,n}=\frac{\left\vert P_{x,x,n}\right\vert }{\left\vert P_{x,n}\right\vert
}%
\]
denote the probability that a random walk of length $n$ starting at $x$ ends
at $x$. Now we prove two general lemmas on vertex transitive graphs that will
be used in the proof of Theorem \ref{kestenaction}.

\begin{lemma}
\label{different}Let $G$ be a $d$-regular vertex transitive graph and let
$n\geq2$ be a positive even integer. Then
\[
\left\vert P_{x,y,n}\right\vert \leq\left\vert P_{x,x,n}\right\vert \leq
d^{2}\left\vert P_{x,x,n-2}\right\vert
\]

\end{lemma}

\noindent\textbf{Proof. } Let $A=dM$ where $M$ is the Markov operator on
$l^{2}(G)$ and let $\delta_{x}$ be the characteristic function of $x$. Then
using the Cauchy-Schwarz inequality we get
\[
\left\vert P_{x,y,n}\right\vert =\left\langle \delta_{x}A^{n/2},\delta
_{y}A^{n/2}\right\rangle \leq\left\langle \delta_{x}A^{n/2},\delta_{x}%
A^{n/2}\right\rangle =\left\vert P_{x,x,n}\right\vert
\]
since $\delta_{x}A^{n}$ is a translate of $\delta_{y}A^{n}$. For the second
inequality, we have
\[
P_{x,x,n}=%
{\displaystyle\bigcup\limits_{y,z\text{ are neighbours of }x}}
P_{y,z,n-2}%
\]
which, using the first inequality and vertex transitivity yields the second
one. $\square$

\begin{lemma}
\label{returningvsrw}Let $G$ be a $d$-regular vertex transitive graph, let
$x\in V(G)$ and let $a$ be a walk of length $l$ starting at $x$. Let $w$ be a
uniform random element of $P_{x,x,n}$. Then the probability
\[
\mathbf{P}\left(  (w_{1},\ldots,w_{l})=a\right)  \geq\frac{1}{d^{2l}}\text{.}%
\]

\end{lemma}

\noindent\textbf{Proof. } The number of walks $w\in P_{x,x,n}$ of length $n$
such that $(w_{1},\ldots,w_{l})=a$ and also $(w_{l+1},\ldots,w_{2l})=a^{-1} $
is exactly $\left\vert P_{x,x,n-2l}\right\vert ,$ which, using Lemma
\ref{different}, is at least
\[
\left\vert P_{x,x,n-2l}\right\vert \geq\frac{1}{d^{2l}}\left\vert
P_{x,x,n}\right\vert
\]
which proves the lemma. $\square$\bigskip

Let $\Gamma$ be a group and let $g$ be a random variable with values in
$\Gamma$. Then $g$ defines the Markov operator $M(g)$ on $l^{2}(\Gamma)$ as
follows. For $f\in l^{2}(\Gamma)$ let
\[
fM(g)=%
{\displaystyle\sum\limits_{h\in\Gamma}}
\boldsymbol{P}(g=h)xh
\]
where $fh$ is the right translate of $f$ by $h$. If $g$ is symmetric, that is,
for all $h\in\Gamma$ we have $P(g=h)=P(g=h^{-1})$, then $M(g)$ is self-adjoint.

\begin{lemma}
\label{modifiedrw}Let $g_{1},g_{2},\ldots,g_{n}$ be independent random
elements of the group $\Gamma$ and let $g=g_{1}g_{2}\cdots g_{n}$. Let $e$
denote the identity element of $\Gamma$. Then
\[
\boldsymbol{P}(g=e)\leq%
{\displaystyle\prod\limits_{i=1}^{n}}
\left\Vert M(g_{i})\right\Vert \text{.}%
\]

\end{lemma}

\noindent\textbf{Proof. } If $g_{1},g_{2}$ are independent, then $M(g_{1}%
g_{2})=M(g_{1})M(g_{2})$. This implies that
\[
M(g)=%
{\displaystyle\prod\limits_{i=1}^{n}}
M(g_{i})\text{.}%
\]
Let $x_{e}$ be the characteristic function of $l^{2}(\Gamma)$. Then, using the
submultiplicativity of norm, we get
\[
\boldsymbol{P}(g=e)=\left\langle x_{e}M(g),x_{e}\right\rangle \leq\left\Vert
x_{e}M(g)\right\Vert \left\Vert x_{e}\right\Vert \leq\left\Vert
M(g)\right\Vert \leq%
{\displaystyle\prod\limits_{i=1}^{n}}
\left\Vert M(g_{i})\right\Vert
\]
as claimed. $\square$

\begin{lemma}
\label{subgroupnorm}Let $\Gamma$ be a group, let $S\subseteq\Gamma$ be a
finite symmetric subset and let $K=\left\langle S\right\rangle $ be the
subgroup of $\Gamma$ generated by $S$. Let $g$ be a uniform random element of
$S$ and let $M$ be the corresponding Markov operator on $l^{2}(\Gamma)$. Then
the norm
\[
\left\Vert M(g)\right\Vert =\rho(\mathrm{Cay}(K,S))\text{.}%
\]

\end{lemma}

\noindent\textbf{Proof. }The space $l^{2}(\Gamma)$ is the orthogonal sum of
countably many isomorphic copies of $l^{2}(K)$ (corresponding to the cosets of
$K$ in $\Gamma$) and $M$ acts diagonally on this space. Hence, the norm equals
the norm of $M$ acting on $l^{2}(K)$, which by definition is the spectral
radius of the corresponding Cayley graph. $\square$\bigskip

A countable group $\Gamma$ is \emph{amenable}, if there is a $\Gamma
$-invariant finitely additive probability measure (called a mean) on all
subsets of $\Gamma$. By Kesten's original theorem (Theorem
\ref{kestenoriginal}), a group $\Gamma$ is amenable, if and only if for all
symmetric probability distributions on $\Gamma$, the corresponding Markov
operator has norm $1$. So a finitely generated group $\Gamma$ is amenable if
and only if $\rho(\mathrm{Cay}(\Gamma,S))=1$ for some (and thus, for every)
finite symmetric generating set $S$ of $\Gamma$.

Amenability is preserved under various group operations, like taking an
ascending union, subgroups or quotient groups or extensions \cite{greenleaf}.
Finite and Abelian (or in general, solvable) groups are amenable and hence the
amenability of $\Gamma$ is preserved under taking a subgroup of finite index.

\section{Invariant Random Subgroups and Probability Measure Preserving Actions
\label{irssection}}

In this section we introduce invariant random subgroups. For a countable group
$\Gamma$ let
\[
\mathrm{Sub}_{\Gamma}=\left\{  H\subseteq\Gamma\mid H\text{ is a subgroup of
}\Gamma\right\}
\]
be the set of subgroups of $\Gamma$, endowed with the product topology
inherited from the space of subsets of $\Gamma$. That is, a sequence of
subgroups $H_{n}$ converges, if for all $g\in\Gamma$, the event $g\in H_{n} $
stabilizes as $n$ tends to infinity. This is called the Chabauty topology
\cite{chaba}. Since $\mathrm{Sub}_{\Gamma}$ is closed in the space of subsets,
$\mathrm{Sub}_{\Gamma}$ is also compact. The group $\Gamma$ acts continuously
on $\mathrm{Sub}_{\Gamma}$ by conjugation.

Assume now that $\Gamma$ is generated by the finite symmetric set $S$. Let
\[
\mathrm{SC}_{\Gamma}(S)=\left\{  \mathrm{Sch}(\Gamma/H,S)\mid H\leq
\Gamma\right\}
\]
be the set of all connected Schreier graphs of $\Gamma$ with respect to $S$.
For $G_{1},G_{2}\in\mathrm{SC}_{\Gamma}(S)$ let $d(G_{1},G_{2})=1/k$ where $k$
is the maximal natural number such that the $k$-balls around the root of
$G_{1}$ and $G_{2}$ are isomorphic as rooted, $S$-labeled graphs. Then $d$ is
a metric on $\mathrm{SC}_{\Gamma}(S)$ and the topology defined by $d$ is
compact. The group $\Gamma$ acts on $\mathrm{SC}_{\Gamma}(S)$ as follows: for
$s\in S$, $\mathrm{Sch}(\Gamma/H,S)s$ is isomorphic to $\mathrm{Sch}%
(\Gamma/H,S)$ but we move the root along the $s$-labeled edge. This extends to
an action of $\Gamma$ on $\mathrm{SC}_{\Gamma}(S)$ by continuous maps. This
can also be expressed as
\[
\mathrm{Sch}(\Gamma/H,S)g=\mathrm{Sch}(\Gamma/H^{g},S)
\]
for $g\in\Gamma$. One can check that the map
\begin{equation}
H\longmapsto\mathrm{Sch}(\Gamma/H,S)\tag{S}%
\end{equation}
is a $\Gamma$-equivariant homeomorphism between $\mathrm{Sub}_{\Gamma}$ and
$\mathrm{SC}_{\Gamma}(S)$ that commutes with the $\Gamma$-action. So, the
spaces $\mathrm{SC}_{\Gamma}(S)$ and $\mathrm{Sub}_{\Gamma}$ are isomorphic as
$\Gamma$-spaces. Note that the inverse of the map (S) can be described as
follows: for $G\in\mathrm{SC}_{\Gamma}(S)$ let $H$ be the set of
$S$-evaluations of returning walks starting at the root in $G$.

\begin{definition}
An \emph{invariant random subgroup} (IRS) of $\Gamma$ is a random subgroup of
$\Gamma$ whose distribution is a $\Gamma$-invariant Borel probability measure
on $\mathrm{Sub}_{\Gamma}$.
\end{definition}

For finitely generated groups, using the above bijection between
$\mathrm{Sub}_{\Gamma}$ and $\mathrm{SC}_{\Gamma}(S)$, every IRS gives rise to
a unimodular random Schreier graph of $\Gamma$, that is, a Borel probability
distribution on $\mathrm{SC}_{\Gamma}(S)$ that is invariant under moving the root.

A natural way to obtain an invariant random subgroup is to take the stabilizer
of a random point in a measure preserving action of $\Gamma$ on a Borel
probability space.

\begin{proposition}
Let $\Gamma$ act on the Borel probability space $(X,\mu)$ by measure
preserving maps. Let $H$ be the stabilizer of a $\mu$-random point in $X$.
Then $H$ is an IRS.
\end{proposition}

\noindent\textbf{Proof. }Let $g\in\Gamma$ be fixed. Then
\[
H^{g}=Stab_{\Gamma}(xg)
\]
where $x$ is $\mu$-random in $X$. But $g$ preserves $\mu$, so $xg$ is also
uniform $\mu$-random in $X$. So the distribution of $H$ and $H^{g}$ are equal.
$\square$\bigskip

Our next proposition says that every IRS actually arises this way.

\begin{proposition}
\label{pmp}Let $H$ be an invariant random subgroup of the finitely generated
group $\Gamma$. Then there exists a measure preserving action of $\Gamma$ on a
Borel probability space $(X,\mu)$ such that $H$ is the stabilizer of a $\mu
$-random point of $X$ in $\Gamma$.
\end{proposition}

\noindent\textbf{Proof. }Fix $S$ to be some finite generating set for $\Gamma
$.\ Let
\[
X=\left\{  (G,f)\mid G\in\mathrm{SC}_{\Gamma}(S),f:V(G)\rightarrow
\lbrack0,1]\right\}
\]
the set of Schreier graphs of $\Gamma$ vertex labeled by elements of the unit
interval. We can endow $X$ with the product topology. The group $\Gamma$
naturally acts on $X$ by moving the root.

Now take our random $H$ and let $G=\mathrm{Sch}(\Gamma/H,S)$ be the Schreier
graph of our random $H$. Label the vertices of this random $G$ using an i.i.d.
with uniform random values in the unit interval $[0,1]$ (according to the
Lebesque measure). This will give a probability measure $\mu$ on $X$ that will
be invariant under the $\Gamma$-action.

Let $B$ be a $\mu$-random element of $X$. Then almost surely, all the labels
of the vertices of $B$ are different. An element $g\in\Gamma$ stabilizes $B$
if and only if when moving the root along a word in $S$ representing $g$ gives
an isomorphic rooted, edge and vertex labeled graph. This implies that $g$
fixes the root, that is, $g\in H$. On the other hand, $H$ trivially stabilizes
$G$. $\square$\bigskip\ 

Note that for non-finitely generated groups, the same proof works: we can take
the random $[0,1]$-labeling simply on the cosets of $H$. We chose to write out
the Schreier graph proof because it is more visual.

Summarizing the above gives us the following.

\begin{proposition}
\label{ekviv}Let $H$ be a random subgroup of the group $\Gamma$, generated by
the finite symmetric set $S$. Then the following are equivalent:\newline%
\qquad1) $H$ is an invariant random subgroup;\newline\qquad2) $\mathrm{Sch}%
(\Gamma/H,S)$ is a unimodular random network; \newline\qquad3) $H$ is the
stabilizer of a random point for some measure preserving action of $\Gamma$.
\end{proposition}

The natural convergence notion for the space of random subgroups is weak
convergence of measures. It turns out that this corresponds to
Benjamini-Schramm convergence \cite{benscha} on the level of Schreier graphs.
For a finite graph $G$ let $\widetilde{G}$ denote the random rooted graph that
we get by assigning the root of $G$ uniformly randomly.

\begin{definition}
Let $G_{n}\in\mathrm{SC}_{\Gamma}(S)$ be a sequence of Schreier graphs and let
$G$ be a random graph in $\mathrm{SC}_{\Gamma}(S)$. We say that $G_{n}$
\emph{Benjamini-Schramm converges to }$G$, if $\widetilde{G_{n}}$ weakly
converges to $G$.
\end{definition}

For a random rooted graph $G$, a finite rooted graph $\alpha$ and $R>0$ let
$P(G,R,\alpha)$ denote the probability that the $R$-ball around the root of
$G$ is isomorphic to $\alpha$. Since the topology on $\mathrm{SC}_{\Gamma}(S)$
is generated by the closed-open sets
\[
U(R,\alpha)=\left\{  G\in\mathrm{SC}_{\Gamma}(S)\mid\text{the }R\text{-ball of
}G\text{ equals }a\right\}
\]
we get that $G_{n}$ Benjamini-Schramm converges to\emph{\ }$G$, if and only
if
\[
\lim_{n\rightarrow\infty}P(G_{n},R,\alpha)=P(G,R,\alpha)
\]
for all $R$ and $\alpha$.

For a group $\Gamma$, a subgroup $H\leq\Gamma$ and $g\in\Gamma$ let
\[
\mathrm{Fix}(\Gamma/H,g)=\left\{  Hx\mid Hxg=Hx\right\}
\]
denote the set of $H$-cosets fixed by $g$.

\begin{lemma}
\label{lekv}Let $\Gamma$ be a group generated by a finite symmetric set $S$
and let $H_{n}\leq\Gamma$ be a sequence of subgroups of finite index with
$\left\vert \Gamma:H_{n}\right\vert \rightarrow\infty$. Then the following are
equivalent: \newline1)\ $\mathrm{Sch}(\Gamma/H_{n},S)$ converges to
$\mathrm{Cay}(\Gamma,S)$ in Benjamini-Schramm convergence; \newline2) $H_{n}$
locally approximates $\Gamma$; \newline3) for every $1\neq g\in\Gamma$, we
have
\[
\lim_{n\rightarrow\infty}\frac{\left\vert \mathrm{Fix}(\Gamma/H_{n}%
,g)\right\vert }{\left\vert \Gamma:H_{n}\right\vert }=0\text{. }%
\]

\end{lemma}

\noindent\textbf{Proof. }Let $G_{n}=\mathrm{Sch}(\Gamma/H_{n},S)$ ($n\geq1 $)
and let $\alpha_{R}$ be the $R$-ball in $\mathrm{Cay}(\Gamma,S)$. Then
\[
P(G_{n},R,\alpha_{R})=\frac{\left\vert \left\{  v\in\mathrm{Sch}(\Gamma
/H_{n},S)\mid B_{R}(v)\cong\alpha_{R}\right\}  \right\vert }{\left\vert
\Gamma:H_{n}\right\vert }%
\]
where $B_{R}(v)$ denotes the ball of radius $R$ centered at $v$. This gives 1)
$\Longleftrightarrow$ 2).

Let $g\in\Gamma$ with $g\neq1$ and let $g=w(S)$ be a word of length $R$
representing $g$ in $S$. Then for every $v\in\mathrm{Fix}(\Gamma/H_{n},g)$,
the path starting at $v$ along the word $w$ is returning, but $w$ does not
return in $\mathrm{Cay}(\Gamma,S)$, so $B_{R}(v)$ is not isomorphic to
$\alpha_{R}$. Hence
\[
\frac{\left\vert \mathrm{Fix}(\Gamma/H_{n},g)\right\vert }{\left\vert
\Gamma:H_{n}\right\vert }\leq1-P(G_{n},R,\alpha_{R})\text{. }%
\]
This yields 1) $\Longrightarrow$ 3).

Vice versa, for a given $R$ and $v\in\mathrm{Sch}(\Gamma/H_{n},S)$, if
$B_{R}(v)$ is not isomorphic to $\alpha_{R}$, then there exists a word $w$ of
length at most $2R$ such that $w$ starting at $v$ is returning but $w(S)\neq
1$. This gives
\[
P(G_{n},R,\alpha_{R})\geq1-\sum\limits_{\left\vert w\right\vert \leq2R\text{,
}w(S)\neq1}\frac{\left\vert \mathrm{Fix}(\Gamma/H_{n},w(S))\right\vert
}{\left\vert \Gamma:H_{n}\right\vert }%
\]
which yields 3) $\Longrightarrow$ 1). $\square$\bigskip\ 

\section{Amenable stabilizers and spectral radius\label{groupsection}}

In this section we prove Theorem \ref{kestenaction}. Note that in the proof,
we are using Kesten's original theorem saying that nonamenable groups have
spectral radius less than $1$, but not Theorem \ref{kestenoriginal}. Hence we
provide an alternate proof for Theorem \ref{kestenoriginal}.

Now we start working towards Theorem \ref{kestenaction}. The next lemma proves
the easy implication. It is essentially the same as Corollary 2 in Kesten's
original paper \cite{kesten1}; since it is only stated there for normal
subgroups, we include a quick proof.

\begin{lemma}
\label{amenablesubgroup}Let $\Gamma$ be a group generated by a finite
symmetric set $S$ and let $H$ be an amenable subgroup of $\Gamma$. Then
\[
\rho(\mathrm{Cay}(\Gamma,S))=\rho(\mathrm{Sch}(\Gamma/H,S)).
\]

\end{lemma}

\noindent\textbf{Proof. } Let $R_{n}$ denote the endpoint of the standard
random $S$-walk on $\Gamma$ starting at the identity $e$. Then $\mathbf{P}%
\left(  R_{n}\in H\right)  $ equals the probability of return for the standard
random $S$-walk on $\Gamma/H$ starting at the coset $H$. Thus, for
$\varepsilon>0$ there exists an integer $n_{0}$ such that for all $n\geq n_{0}
$ we have
\[
q_{2n}=\mathbf{P}\left(  R_{2n}\in H\right)  \geq\left(  (1-\varepsilon
)\rho\right)  ^{2n}%
\]
where $\rho=\rho(\mathrm{Sch}(\Gamma/H,S))$. Now let $n>n_{0}$ to be chosen
later and for $h\in H$ let
\[
Q(h)=\frac{\mathbf{P}\left(  R_{2n}=h\right)  }{q_{2n}}%
\]
Then $Q$ is a finite symmetric probability distribution on $H$ and since $H$
is amenable, using the original Kesten's theorem, the random walk $R^{\prime}
$ on $H$ with respect to $Q$ has spectral radius $1$. In particular, there
exists an integer $m$ such that
\[
\mathbf{P}\left(  R_{2m}^{\prime}=e\right)  \geq(1-\varepsilon)^{2m}\text{.}%
\]
Now let us consider the probability, that the walk $R$ of length $4nm$ returns
to $H$ in every $2n$-th step and the product of the segments is $e$. This
gives
\[
\mathbf{P}\left(  R_{4nm}=e\right)  \geq q_{2n}^{2m}\mathbf{P}\left(
R_{2m}^{\prime}=e\right)  \geq\rho^{4nm}(1-\varepsilon)^{2m(2n+1)}%
\]
Taking $4nm$-t roots, and using that $n$ is arbitrarily large, we get
\[
\rho(\mathrm{Cay}(\Gamma,S))\geq(1-\varepsilon)\rho
\]
Hence, the lemma holds. $\square$\bigskip

Let $\Gamma$ be a group and let $S\subseteq\Gamma$ be a finite symmetric
multiset. Let $H$ be a subgroup of $\Gamma$. Let $S^{n}$ denote the set of $n
$-tuples in $S$ and for $a\in S^{n}$ let $[a]=a_{1}a_{2}\cdots a_{n}\in\Gamma$
be the product of elements in the tuple. Let
\[
A_{H}(S,n)=\left\{  a\in S^{n}\mid\lbrack a]\in H\right\}  \text{.}%
\]
Let $A(S,n)=A_{1}(S,n)$. The set $A(S,n)$ can be identified with returning
walks of length $n$ in $\mathrm{Cay}(\Gamma,S)$.

\begin{lemma}
\label{triv2}Let $(a_{0},\ldots,a_{n-1})$ be a uniform random element of
$A(S,n)$. Then for any $1\leq k$ and $0\leq t\leq n-1$, the distribution of
the segment $(a_{t},\ldots,a_{t+k-1})\in S^{k}$, where the indices are
understood modulo $n$, is independent of $l$.
\end{lemma}

\noindent\textbf{Proof. } The set $A(S,n)$ is invariant under cyclic
rotations. $\square$\ \bigskip

\begin{lemma}
\label{triv1}Let $k>0$ be an integer and let $(a_{0},\ldots,a_{n-1})$ be a
uniform random element of $A(S,n)$ with $n>2k$. Then for any $b\in S^{k}$, the
probability
\[
\mathbf{P}\left(  (a_{0},a_{1},\ldots,a_{k-1})=b\right)  \geq\left\vert
S\right\vert ^{-2k}\text{.}%
\]

\end{lemma}

\noindent\textbf{Proof. } Let $b=(b_{0},b_{1},\ldots,b_{k-1})$. Then for any
$(c_{0},\ldots,c_{n-2k-1})\in A(S,n-2k)$, we have
\[
(b_{0},\ldots,b_{k-1},b_{k-1}^{-1},\ldots,b_{0}^{-1},c_{0},c_{1}%
,\ldots,c_{n-2k-1})\in A(S,n)
\]
which gives
\[
\mathbf{P}\left(  (a_{0},a_{1},\ldots,a_{k-1})=b\right)  \geq\frac{\left\vert
A(S,n-2k)\right\vert }{\left\vert A(S,n)\right\vert }=\frac{p_{n-2k}}{p_{n}%
}\left\vert S\right\vert ^{-2k}\geq\left\vert S\right\vert ^{-2k}%
\]
where $p_{n}$ denotes the probability of return on $\mathrm{Cay}(\Gamma,S)$ in
$n$ steps. In the last inequality we are using Lemma \ref{different}.
$\square$\bigskip

For $a\in S^{n}$ we shall look at the event
\[
Ha_{0}\cdots a_{t-1}=Ha_{0}\cdots a_{t}%
\]
as right cosets of $H$. Equivalently, the walk corresponding to $a$ in
$\mathrm{Sch}(\Gamma/H,S)$ is at the same coset at time $t-1$ and $t$. Yet
another way of writing this is that $a_{t}\in H^{a_{0}\cdots a_{t-1}}$.

For $a\in S^{n}$ let the index set
\[
I(S,H,a)=\left\{  t\in\{0,\ldots,n-1\}\mid a_{t}\in H^{a_{0}\cdots a_{t-1}%
}\right\}  \text{.}%
\]

\begin{definition}
We call $a,b\in S^{n}$ $H$-\emph{equivalent} if
\[
Ha_{0}\cdots a_{t}=Hb_{0}\cdots b_{t}\text{ }(0\leq t\leq n-1)\text{ and
}a_{t}=b_{t}\text{ (}t\notin I(S,H,a)\text{)}%
\]
Let $C(a)$ denote the $H$-equivalence class of $a$.
\end{definition}

Note that $H$-equivalent sequences have the same index set.

\begin{lemma}
\label{inverz}Let $\Gamma$ be a group and let $S\subseteq\Gamma$ be a finite
symmetric multiset. Let $H$ be a subgroup of $\Gamma$ and let $n>0$ be an
integer. Then
\[
\left\vert A_{H}(S,n)\right\vert \geq%
{\displaystyle\sum\limits_{a\in A(S,n)}}
\left[  p(S,H,a)\right]  ^{-1}%
\]
where
\[
p(S,H,a)=%
{\displaystyle\prod\limits_{t\in I(S,H,a)}}
\rho(\mathrm{Cay}(H^{a_{0}\cdots a_{t-1}},S\cap H^{a_{0}\cdots a_{t-1}%
}))\text{.}%
\]

\end{lemma}

\noindent\textbf{Proof. } For $0\leq t\leq n-1$ let
\[
X_{t}=S\cap H^{a_{0}\cdots a_{t-1}}\text{.}%
\]
Then $X_{t}$ is a symmetric multiset. The definition of $H$-equivalence
implies that $C(a)\subseteq A_{H}(S,n)$ for all $a\in A_{H}(S,n)$ and thus,
$A_{H}(S,n)$ is the disjoint union of its $H$-equivalence classes. A uniform
random element $x$ of $C(a)$ is of the form $\left(  x_{1},x_{2},\ldots
,x_{n}\right)  $ where $x_{t}=a_{t}$ is fixed for $t\notin I(S,H,a)$ and
$x_{t}$ is a uniform random element of $X_{t}$ for $t\in I(S,H,a)$. Hence,
using Lemma \ref{modifiedrw} and Lemma \ref{subgroupnorm}, for a uniform
random element $x$ of $C(a)$ we have
\[
\boldsymbol{P}\left(  x\in A(S,n)\right)  \leq%
{\displaystyle\prod\limits_{t\in I(S,H,a)}}
\rho(\mathrm{Cay}(H^{a_{0}\cdots a_{t-1}},S\cap H^{a_{0}\cdots a_{t-1}%
}))=p(S,H,a)\text{.}%
\]

We can estimate the size of $A_{H}(S,n)$ by the sum of the sizes of
equivalence classes of tuples in $A(S,n)$. We count $C(a)$ exactly $\left\vert
C(a)\right\vert \boldsymbol{P}\left(  x\in A(S,n)\right)  $ times, which
gives
\[
\left\vert A_{H}(S,n)\right\vert \geq%
{\displaystyle\sum\limits_{a\in A(S,n)}}
\left[  \boldsymbol{P}(x\in A(S,n))\right]  ^{-1}\geq%
{\displaystyle\sum\limits_{a\in A(S,n)}}
\left[  p(S,H,a)\right]  ^{-1}\text{.}%
\]
The lemma holds. $\square$\bigskip

Let us denote
\[
\lbrack S^{k}]=\left\{  a_{0}\cdots a_{k-1}\mid a\in S^{k}\right\}
\]
where we look at $[S^{k}]$ as a multiset in $\Gamma$. Let $\{S^{k}\}$ denote
the set of elements in $[S^{k}]$.

\begin{lemma}
\label{findsubgroup}Let $\Gamma$ be a group generated by the finite symmetric
set $S$ and let $H$ be an invariant random subgroup of $\Gamma$ that is
nonamenable with positive probability. Then there exists an integer $k>0$ and
$p,r>0$ such that
\[
\rho(\mathrm{Cay}(H,[S^{k}]\cap H))<1-r
\]
with probability at least $p$.
\end{lemma}

\noindent\textbf{Proof. } By adding $ss^{-1}$ to words, we see that for
$k\geq0$ we have $\{S^{k}\}\subseteq\{S^{k+2}\}$. Since $S$ generates $\Gamma
$, the subgroup
\[
N=%
{\displaystyle\bigcup\limits_{k=0}^{\infty}}
\{S^{2k}\}
\]
has index at most $2$ in $\Gamma$. Let $H_{2k}$ be the subgroup generated by
$\{S^{2k}\}\cap H$. The group $H_{2k}$ is a well-defined subgroup of the
random subgroup $H$ of $\Gamma$, and the union $%
{\displaystyle\bigcup}
H_{2k}$ has index at most $2$ in $H$. Hence, there exists $k$ such that
$H_{2k}$ is non-amenable with positive probability, using that the ascending
union of amenable groups is amenable and that amenability does not change when
passing to a finite index subgroup. By Kesten's theorem on amenability
\cite{kesten1}, $H_{2k}$ is nonamenable if and only if $\rho(\mathrm{Cay}%
(H_{2k},[S^{2k}]\cap H))<1$. Since $\rho(\mathrm{Cay}(H,[S^{2k}]\cap H))$ is a
measurable function of $H$, there exists $p,r>0$ such that $\rho
(\mathrm{Cay}(H_{2k},[S^{2k}]\cap H))<1-r$ with probability at least $p$.
$\square$\bigskip

We are ready to prove the main theorem of this section after a trivial lemma.

\begin{lemma}
\label{trivkancs}Let $X$ be a random variable with $0\leq X\leq R$. Then
\[
\boldsymbol{P}\left(  X\geq\frac{\boldsymbol{E}\left[  X\right]  }{2}\right)
\geq\frac{\boldsymbol{E}\left[  X\right]  /R}{2-\boldsymbol{E}\left[
X\right]  /R}\text{.}%
\]

\end{lemma}

\bigskip

\noindent\textbf{Proof of Theorem \ref{kestenaction}. }If $H$ is amenable
a.s., then by Lemma \ref{amenablesubgroup}, we have
\[
\rho(\mathrm{Sch}(\Gamma/H,S))=\rho(\mathrm{Cay}(\Gamma,S))
\]
a.s. This part of the theorem does not require invariance or random subgroups.

Assume $H$ is nonamenable with positive probability. Then by Lemma
\ref{findsubgroup} there exists an integer $k>0$ and $p,r>0$ such that
\[
\rho(\mathrm{Cay}(H,[S^{k}]\cap H))<1-r
\]
with probability at least $p$.

Let $T=[S^{k}]$. Fix a positive integer $n$. Let
\[
J(T,H,a)=\left\{  t\in I(T,H,a)\mid\rho(\mathrm{Cay}(H^{a_{0}\cdots a_{t-1}%
},T\cap H^{a_{0}\cdots a_{t-1}}))<1-r\right\}  \text{.}%
\]
Then we have
\[
p(T,H,a)=%
{\displaystyle\prod\limits_{t\in I(a)}}
\rho(\mathrm{Cay}(H^{a_{0}\cdots a_{t-1}},T\cap H^{a_{0}\cdots a_{t-1}%
}))<(1-r)^{\left\vert J(T,H,a)\right\vert }%
\]
and hence, by Lemma \ref{inverz}, for any $H\leq\Gamma$ we have
\begin{equation}
\left\vert A_{H}(T,n)\right\vert \geq%
{\displaystyle\sum\limits_{a\in A(T,n)}}
p^{-1}(T,H,a)>%
{\displaystyle\sum\limits_{a\in A(T,n)}}
(1-r)^{-\left\vert J(T,H,a)\right\vert }\text{.}\tag{A}%
\end{equation}

For any fixed element $a\in A(T,n)$ and $0\leq t\leq n-1$, over the random
subgroup $H$, the probability
\begin{align*}
\boldsymbol{P}\left(  t\in J(T,H,a)\right)   &  =\\
&  =\boldsymbol{P}\left(  a_{t}\in H^{a_{0}\cdots a_{t-1}}\text{, }%
\rho(\mathrm{Cay}(H^{a_{0}\cdots a_{t-1}},T\cap H^{a_{0}\cdots a_{t-1}%
}))<1-r\right)  =\\
&  =\boldsymbol{P}\left(  a_{t}\in H\text{ and }\rho(\mathrm{Cay}(H,T\cap
H))<1-r\right)
\end{align*}
by the invariance of $H$. Thus, the expected value of $\left\vert
J(T,H,a)\right\vert $ over the random subgroup $H$ equals
\begin{equation}
\boldsymbol{E}\left[  \left\vert J(T,H,a)\right\vert \right]  =%
{\displaystyle\sum\limits_{t=0}^{n-1}}
\boldsymbol{P}\left(  a_{t}\in H\text{ and }\rho(\mathrm{Cay}(H,T\cap
H))<1-r\right)  \text{.}\tag{B}%
\end{equation}

For any fixed subgroup $H$ such that $\rho(\mathrm{Cay}(H,T\cap H))<1-r$ we
have $T\cap H\neq\varnothing$. In particular, for these $H$, using Lemma
\ref{triv2} and Lemma \ref{triv1}, we have
\[
\frac{1}{\left\vert A(T,n)\right\vert }%
{\displaystyle\sum\limits_{a\in A(T,n)}}
\boldsymbol{1}\left(  a_{t}\in H\right)  \geq\left\vert T\right\vert ^{-2}%
\]
for all $0\leq t\leq n-1$. Since $\boldsymbol{P}\left(  \rho(\mathrm{Cay}%
(H,T\cap H))<1-r\right)  >p$, we get
\[
\frac{1}{\left\vert A(T,n)\right\vert }%
{\displaystyle\sum\limits_{a\in A(T,n)}}
\boldsymbol{P}\left(  a_{t}\in H\text{ and }\rho(\mathrm{Cay}(H,T\cap
H))<1-r\right)  \geq p\left\vert T\right\vert ^{-2}%
\]
and so by summing with respect to $t$ and using (B) we yield
\[
\boldsymbol{E}\left[  \frac{1}{\left\vert A(T,n)\right\vert }%
{\displaystyle\sum\limits_{a\in A(T,n)}}
\left\vert J(T,H,a)\right\vert \right]  \geq p\left\vert T\right\vert
^{-2}n\text{.}%
\]

Using
\[
\frac{1}{\left\vert A(T,n)\right\vert }%
{\displaystyle\sum\limits_{a\in A(T,n)}}
\left\vert J(T,H,a)\right\vert \leq n
\]
and Lemma \ref{trivkancs} we get%

\begin{equation}
\boldsymbol{P}\left(  \frac{1}{\left\vert A(T,n)\right\vert }%
{\displaystyle\sum\limits_{a\in A(T,n)}}
\left\vert J(T,H,a)\right\vert \geq\frac{1}{2}p\left\vert T\right\vert
^{-2}n\right)  \geq\frac{p\left\vert T\right\vert ^{-2}}{2-p\left\vert
T\right\vert ^{-2}}\text{.}\tag{C}%
\end{equation}

Let $H$ be a subgroup satisfying (C). Using the inequality of arithmetic and
geometric means and (C) we get
\[
\frac{1}{\left\vert A(T,n)\right\vert }%
{\displaystyle\sum\limits_{a\in A(T,n)}}
(1-r)^{-\left\vert J(T,H,a)\right\vert }\geq(1-r)^{-\frac{1}{\left\vert
A(T,n)\right\vert }%
{\displaystyle\sum\limits_{a\in A(T,n)}}
\left\vert J(T,H,a)\right\vert }\geq(1-r)^{-\frac{1}{2}p\left\vert
T\right\vert ^{-2}n}%
\]
which by (A) gives
\[
\left\vert A_{H}(T,n)\right\vert >%
{\displaystyle\sum\limits_{a\in A(T,n)}}
(1-r)^{-\left\vert J(T,H,a)\right\vert }\geq\left\vert A(T,n)\right\vert
(1-r)^{-\frac{1}{2}p\left\vert T\right\vert ^{-2}n}.
\]
Using $\left\vert A_{H}(T,n)\right\vert =\left\vert A_{H}(S,kn)\right\vert $,
$\left\vert A(T,n)\right\vert =\left\vert A(S,kn)\right\vert $, dividing by
$d^{kn}$ and taking the $nk$-th roots of both sides, we get the following.

For every integer $n>1$, with probability at least $p\left\vert T\right\vert
^{-2}/(2-p\left\vert T\right\vert ^{-2})$, our random $H$ satisfies%

\[
\left(  p_{nk}(\mathrm{Sch}(\Gamma/H,S))\right)  ^{1/nk}>\left(
p_{nk}(\mathrm{Cay}(\Gamma,S))\right)  ^{1/nk}(1-r)^{-\frac{1}{2k}p\left\vert
T\right\vert ^{-2}}%
\]
Since
\[
\lim_{n\rightarrow\infty}\left(  p_{2nk}(\mathrm{Cay}(\Gamma,S))\right)
^{1/2nk}=\rho(\mathrm{Cay}(\Gamma,S)
\]
we get that there exists $n$ for which
\[
\rho(\mathrm{Sch}(\Gamma/H,S))\geq\left(  p_{nk}(\mathrm{Sch}(\Gamma
/H,S))\right)  ^{1/nk}>\rho(\mathrm{Cay}(\Gamma,S)
\]
with probability at least $p\left\vert T\right\vert ^{-2}/(2-p\left\vert
T\right\vert ^{-2})$. The theorem holds. $\square$

\section{Local and spectral approximation for linear groups}

In this section we prove Proposition \ref{abo}, Theorem \ref{farber} and
Theorem \ref{raman}. We start with an easy lemma on free groups.

\begin{lemma}
\label{free}Let $H$ be an amenable invariant random subgroup of the countable
non-Abelian free group $F$. Then $H=1$ a.s.
\end{lemma}

\noindent\textbf{Proof. }By the Nielsen-Schreier theorem, all subgroups of $F$
are free. Since non-Abelian free groups are non-amenable, we get that $H$ is
cylic a.s. The group $F$ has only countably many cyclic subgroups, and for
every infinite cyclic subgroup $C$ of $F$, the conjugacy class of $C$ in $F$
is infinite. Hence, every invariant measure supported on cyclic subgroups is
supported on the trivial subgroup. $\square$\ \bigskip

The following lemma is a variant on the Alon-Boppana theorem \cite{alon}.

\begin{lemma}
\label{albo}Let $\Gamma$ be an infinite group generated by the finite
symmetric set $S$ and let $H_{n}\leq\Gamma$ be a sequence of subgroups of
finite index with $\left\vert \Gamma:H_{n}\right\vert \rightarrow\infty$ such
that $H_{n}\rightarrow H$ in Chabauty topology. Then
\[
\lim\inf\rho_{0}(\mathrm{Sch}(\Gamma/H_{n},S))\geq\rho(\mathrm{Sch}%
(\Gamma/H,S))\text{.}%
\]

\end{lemma}

\noindent\textbf{Proof. }By Chabauty convergence, $\mathrm{Sch}(\Gamma
/H_{n},S)$ converges to $\mathrm{Sch}(\Gamma/H,S)$ in rooted convergence (see
Section \ref{irssection}). Since $\left\vert \mathrm{Sch}(\Gamma
/H_{n},S)\right\vert \rightarrow\infty$, we get that $H$ has infinite index in
$\Gamma$. By the definition of the spectral radius as a norm, for any
$\varepsilon>0$ there exists a finitely supported function $f\in l^{2}%
(\Gamma/H)$ with zero sum such that $\left\langle f,f\right\rangle =1$ and
\[
\left\langle fM,f\right\rangle >\rho(\mathrm{Sch}(\Gamma/H,S))-\varepsilon
\text{.}%
\]
Here $M$ denotes the Markov operator on $l^{2}(\Gamma/H)$.

Let $R>0$ such that the ball $B$ of radius $R$ in $\mathrm{Sch}(\Gamma/H,S)$
centered at the root contains the support of $f$. By rooted convergence, for
every large enough $n$, the $R$-ball $B^{\prime}$ in $\mathrm{Sch}%
(\Gamma/H,S)$ centered at the root is isomorphic to $B$. Let $f^{\prime}\in
l^{2}(\Gamma/H_{n})$ be the same as $f$ on $B^{\prime}$ (using the isomorphism
between $B$ and $B^{\prime}$) and $0$ outside. Then $f^{\prime}\in l_{0}%
^{2}(\Gamma/H_{n})$, $\left\langle f^{\prime},f^{\prime}\right\rangle =1$ and
\[
\rho_{0}(\mathrm{Sch}(\Gamma/H_{n},S))\geq\left\langle f^{\prime}M,f^{\prime
}\right\rangle =\left\langle fM,f\right\rangle >\rho(\mathrm{Sch}%
(\Gamma/H,S))-\varepsilon\text{.}%
\]
This gives
\[
\lim\inf\rho_{0}(\mathrm{Sch}(\Gamma/H_{n},S))\geq\rho(\mathrm{Sch}%
(\Gamma/H,S))-\varepsilon
\]
which proves the lemma. $\square$\bigskip

\noindent\textbf{Proof of Proposition \ref{abo}. }By the compactness of the
Chabauty topology, any subsequence of $H_{n}$ has a Chabauty convergent
subsequence $K_{n}$ with limit $K$. Applying Lemma \ref{albo} gives
\[
\lim\inf\rho_{0}(\mathrm{Sch}(\Gamma/K_{n},S))\geq\rho(\mathrm{Sch}%
(\Gamma/K,S))\geq\rho(\mathrm{Cay}(\Gamma,S))
\]
using $\rho(\mathrm{Sch}(\Gamma/K,S))\geq\rho(\mathrm{Cay}(\Gamma,S))$. This
proves the proposition. $\square$\bigskip

We are ready to prove Theorem \ref{farber}. In the proof, we are using the
following result from \cite{abgv}. \bigskip

\begin{theorem}
\label{idez}Let $\Gamma$ be a finitely generated linear group and let $H$ be
an IRS of $\Gamma$ that is amenable a.s. Then $H\leq A(\Gamma)$ a.s., where
$A(\Gamma)$ denotes the amenable radical of $\Gamma$. 
\end{theorem}

\noindent\textbf{Proof of Theorem \ref{farber}. }Assume by contradiction that
$H_{n}$ does not approximate $\Gamma$ locally. Then by Lemma \ref{lekv} there
exists $1\neq g\in\Gamma$, $c>0$ and a subsequence $K_{n}$ of $H_{n}$ such
that we have
\[
\left\vert \mathrm{Fix}(\Gamma/K_{n},g)\right\vert >c\left\vert \Gamma
:K_{n}\right\vert \text{ \ }(n\geq1)\text{ }%
\]
where
\[
\mathrm{Fix}(\Gamma/K_{n},g)=\left\{  K_{n}x\mid K_{n}xg=K_{n}x\right\}
\]
is the set of $H$-cosets fixed by $g$. Note that the above condition is
equivalent to $g\in H^{x}$. Let $M_{n}$ denote a uniform random conjugate of
$K_{n}$. Then the invariant random subgroup $M_{n}$ equals the stabilizer of a
uniform random vertex in $\mathrm{Sch}(\Gamma/K_{n},S)$, so we have $g\in
M_{n}$ with probability at least $c$ ($n\geq1$).

By passing to a subsequence of $K_{n}$, we can assume that $M_{n}$ is weakly
convergent. Let $M$ be the limit of $M_{n}$. By weak convergence, we have
$g\in M_{n}$ with probability at least $c$.

Let $L$ be in the support of $M$. Then by weak convergence, there exists
$L_{n}$ in the support of $M_{n}$ with%
\[
\lim_{n\rightarrow\infty}L_{n}=L\text{ }%
\]
in the Chabauty topology. Using Lemma \ref{albo} we get
\[
\rho(\mathrm{Sch}(\Gamma/L,S))\leq\lim\inf\rho_{0}(\mathrm{Sch}(\Gamma
/L_{n},S))\leq\rho(\mathrm{Cay}(\Gamma,S))
\]
by the assumption of the theorem, since $L_{n}$ is a conjugate of $K_{n}$ and
hence $\mathrm{Sch}(\Gamma/L_{n},S)$ is isomorphic to $\mathrm{Sch}%
(\Gamma/K_{n},S)$. So, $M$ is an IRS of $\Gamma$ satisfying $\rho
(\mathrm{Sch}(\Gamma/M,S))=\rho(\mathrm{Cay}(\Gamma,S))$ a.s., and so by
Theorem \ref{kestenaction}, $M$ is amenable a.s.

Since $\Gamma$ is linear, by Theorem \ref{idez}, every amenable IRS of
$\Gamma$ lies in the amenable radical of $\Gamma$. By the assumption of the
theorem, this is trivial, so $M=1$ a.s. But $g\in M$ with probability at least
$c$, a contradiction. $\square$\bigskip

Now Theorem \ref{raman} follows fast. \bigskip

\noindent\textbf{Proof of Theorem \ref{raman}. }Let $\phi:F_{S}\rightarrow
\Gamma$ be the canonical projection defined by $S$ and let $K_{n}=\phi
^{-1}(H_{n})$. Then $\mathrm{Sch}(\Gamma/H_{n},S)$ is isomorphic to
$\mathrm{Sch}(F_{S}/K_{n},S)$. Since free groups are linear, we can apply
Theorem \ref{farber} to $K_{n}$ and get that $K_{n}$ locally approximates
$F_{S}$. In particular, $\mathrm{Sch}(\Gamma/H_{n},S)$ has essentially large
girth. Let $M_{n}$ denote a uniform random conjugate of $K_{n}$. Then $M_{n}$
converges to the trivial group, but $M_{n}$ contains the kernel of $\phi$ a.s.
($n\geq1$). This implies that $\phi$ is an isomorphism and hence $\Gamma$ is a
free group freely generated by $S$. $\square$\bigskip

Note that in the proof of Theorem \ref{raman} one can use the easier Lemma
\ref{free} instead of Theorem \ref{idez}, which makes it independent from
\cite{abgv}. \bigskip

\noindent\textbf{Acknowledgments.} M.A.\ has been supported by the grant
IEF-235545. Y.G.\ was partially supported by ISF grant 441/11 B.V.\ was
supported by the NSERC Discovery Accelerator Grant and the Canada Research
Chair program.

\noindent\textsc{Mikl\'os Ab\'ert.} Alfr\'ed R\'enyi Institute of Mathematics.
Re\'altanoda utca 13-15, H-1053, Budapest, Hungary. \texttt{abert\@@renyi.hu}
\bigskip

\noindent\textsc{Yair Glasner.} Department of Mathematics. Ben-Gurion
University of the Negev. P.O.B. 653, Be'er Sheva 84105, Israel.
\texttt{yairgl\@@math.bgu.ac.il}\bigskip

\noindent\textsc{B\'{a}lint Vir\'{a}g.} Departments of Mathematics and
Statistics. University of Toronto, M5S 2E4 Canada.
\texttt{balint\@@math.toronto.edu}

\end{document}